\documentclass[12pt,reqno]{article}

\usepackage[top=2.5cm,bottom=2.5cm,left=2.5cm,right=2.5cm]{geometry}
\usepackage{graphicx}
\usepackage{latexsym}
\def\Dj{\hbox{D\kern-.73em\raise.30ex\hbox{-}
\raise-.30ex\hbox{}}}
\def\dj{\hbox{d\kern-.33em\raise.80ex\hbox{-}
\raise-.80ex\hbox{\kern-.40em}}}
\usepackage{amssymb}
\usepackage{amscd}
\usepackage[arrow,matrix]{xy}
\usepackage{graphicx}
\usepackage{xcolor}
\usepackage{times}
\usepackage{mathptmx}
\usepackage[]{subfigure}
\RequirePackage{hyperref}
\RequirePackage{bookmark}
 \usepackage{tabu}
 \usepackage{cite}
 \usepackage{float}

\usepackage{amsfonts}

\usepackage{amscd,amsmath,amsthm,amssymb}

\theoremstyle{plain}
\newtheorem{theorem}{Theorem}[section]

\newtheorem{Property}{Property}[section]
\newtheorem{proposition}[Property]{Proposition}
\newtheorem{corollary}[Property]{Corollary}
\theoremstyle{definition}

\newtheorem{example}[Property]{Example}

\newtheorem{lemma}[Property]{Lemma}

\begin{document}

\title{\textbf{Elliptical trajectories of a point on the elliptical 2-sphere }}
\author {Zehra \"{O}zdemir$^{1}${}\thanks{{
 E--mail: \texttt{zehra.ozdemir@amasya.edu.tr } (Z. \"{O}zdemir)}},\texttt{ }  Fatma Ate\c{s}$^{2}$\  \\
 %EndAname
\begin{small}{$^{1}$Department of Mathematics, Arts and Science Faculty, Amasya University, 05189 Amasya, Turkey}\end{small}\\
\begin{small}{$^{2}$Department of Mathematics-Computer, Arts and Science Faculty, Necmettin Erbakan University}\end{small}, \begin{small}{Konya,
Turkey}\end{small}}
\date{}
\maketitle
%%%%%%%%%%%%%%%%%%%%%%%%%%%%%%
\begin{abstract} 
The focus of this work is to analyze the trajectories of a point on the
ellipsoid $\mathbb{S}_{a_{1},a_{2},a_{3}}^{2}$ while it is under the
influence of a Killing vector field $K$. For this purpose, we introduce the generalized Darboux frame and the variational vector fields of $\mathbb{S}_{a_{1},a_{2},a_{3}}^{2}$. Then, we determine the Killing equations in terms of the Darboux frame invariants along an ellipsoidal curve. The Killing equations make it possible for us to interpret the magnetic trajectory of a point on the ellipsoid $\mathbb{S}_{a_{1},a_{2},a_{3}}^{2}$. Then, we determine two special trajectories using the variational method. The first one is magnetic curves that are the trajectories produced by the Killing magnetic field $K$ are satisfied the following Lorentz force equation $F_{L} (t)=K\times _{E}t=\nabla _{T}t$, where $\times _{E}$ is elliptical cross product and $\nabla $ is the Levi-Civita connection of the ellipsoid $\mathbb{S}_{a_{1},a_{2},a_{3}}^{2}$. The second one is generalized magnetic helices that are trajectories described by the trajectory of a point on a great ellipse of the ellipsoid rolling without slipping on a fixed ellipse of the ellipsoid using the elliptical motion on the $\mathbb{S}_{a_{1},a_{2},a_{3}}^{2}$. Furthermore, we give various examples and visualized them with the program Mathematica.
\end{abstract}
%%%%%%%%%%%%
\begin{small} {\textbf{MSC:} 53Z04, 53B50, 37C10, 14H45, 14H50, 35A15, 70E17.}
\end{small}\\
\begin{small} {\textbf{Keywords:} Applications to physics, Vector fields, Magnetic flows, Ordinary differential equations, Special curves, Variational methods, Motion of a rigid body with a fixed point.} 
\end{small}\\

\section{Introduction}

The motion of a rigid body around an axis in three dimensions is a linear transformation that can be expressed with an orthonormal matrix. This matrix called a rotation matrix. The rotation matrices are used extensively for computations in geometry, kinematics, physics, computer graphics, animations, and optimization problems involving the estimation of rigid body transformations. With the help of these motions, some special curves are obtained. For example; spherical helices are trajectories of a point on a great circle of a sphere while rolling without slipping on a fixed circle of the sphere. On the other hand, the elliptical rotation is defined as the motion of a point on an ellipse through some angle about a given vector. In \cite{Ozd1,Ozd2}, the author generated the elliptical and  hyperbolical rotation matrices using the Rodrigues, Cayley and Householder methods on the ellipsoid (elliptical 2-sphere) and hyperboloid.  

The purpose of this paper is to determine the trajectories of a point on the ellipsoid by using the elliptical structures. For this, we analyze the magnetic trajectories of a point on an ellipsoid using the elliptical motion. So, we define Darboux frame equations associated with the elliptical inner product and the elliptical vector product. Then, we compute the variations of the elliptical Darboux frame equations. After that, we obtain the charged particle trajectories on the elliptical 2-sphere by using these variations. Finally, we give various examples and illustrate their images on the elliptical 2-sphere.

The paper is organized as follows: The first section is reserved for introduction. Some basic concepts are presented in the second section. The variational formulations of the curves on the $\mathbb{S}_{a_{1},a_{2},a_{3}}^{2}$ are investigated in the third section. Also, in this section, the relationships between the Killing equations and the variation formulas of Darboux frame curvatures along an ellipsoidal curve are examined. Then using the Killing equations the magnetic trajectories of a point on the $\mathbb{S}_{a_{1},a_{2},a_{3}}^{2}$ are obtained and some examples are given in the fourth section. In the fifth section, helical trajectories of the point on the elliptical 2-sphere are investigated by using the variational method and elliptical motion. In this section also some elliptical trajectories such as satellite and cycloid are investigated via elliptical motions. Besides, some related examples are plotted. Finally, the physical and geometrical interpretations of the results we have investigated are presented in the last section.

\section{Basic Concepts}

In this section, we present brief information of the structures on the $%
\mathbb{S}_{a_{1},a_{2},a_{3}}^{2}$ to describe the basic background.

The real vector space $\mathbb{R}^{3}$ furnished with the elliptical inner
product
\begin{equation*}
B:\mathbb{R}^{3}\times \mathbb{R}^{3}\rightarrow \mathbb{R};\text{ }%
B(x,y)=a_{1}x_{1}y_{1}+a_{2}x_{2}y_{2}+a_{3}x_{3}y_{3}.
\end{equation*}%
is represented by $\mathbb{R}_{a_{1},a_{2},a_{3}}^{3}.$ Where $%
x=(x_{1},x_{2},x_{3}),$ $y=(y_{1},y_{2},y_{3})\in \mathbb{R}^{3}$ and $%
a_{1},a_{2},a_{3}\in \mathbb{R}^{+}.$ The norm of a vector $x$, associated with
the scalar product $B$, is defined as follows
\begin{equation*}
\left\Vert x\right\Vert _{B}=\sqrt{B(x,x)}.
\end{equation*}%
Two vectors $x$ and $y$ are called elliptically orthogonal vectors if $%
B(x,y)=0.$ In addition, if $x$ is a elliptically orthonormal vector then $%
B(x,x)=1$. The cosine of the angle between two vectors $x$ and $y$ is
defined as%
\begin{equation*}
\cos \theta =\frac{B(x,y)}{\left\Vert x\right\Vert _{B}\left\Vert
y\right\Vert _{B}},
\end{equation*}%
where $\theta $ is compatible with the parameters of the angular parametric
equations of ellipse or elliptical 2-sphere \cite{Ozd1}.

The elliptical cross product of two vector fields in $\mathbb{R}%
_{a_{1},a_{2},a_{3}}^{3}$ is given by

\begin{equation}
X\times _{E}Y=\Delta \left\vert
\begin{array}{ccc}
\frac{e_{1}}{a_{1}} & \frac{e_{2}}{a_{2}} & \frac{e_{3}}{a_{3}} \\
x_{1} & x_{2} & x_{3} \\
y_{1} & y_{2} & y_{3}%
\end{array}%
\right\vert,
\end{equation}%
where $\Delta =\sqrt{a_{1}a_{2}a_{3}},$ $a_{1},a_{2},a_{3}\in \mathbb{R}^{+}$%
\cite{Ozd1}.

A surface of center $p_{0}$ and radius $r$ in $\mathbb{R}%
_{a_{1},a_{2},a_{3}}^{3}$ is defined by the form%
\begin{equation*}
S_{a_{1},a_{2},a_{3}}^{2}(p_{0};r)=\{p\in \mathbb{R}%
^{3}:B(p-p_{0},p-p_{0})=r^{2}\}.
\end{equation*}%
When $p_{0}$ is the origin and $r=1$ the surface called as the elliptical 2
sphere and denoted by%
\begin{equation*}
S_{a_{1},a_{2},a_{3}}^{2}=\{(x,y,z)\in \mathbb{R}%
^{3}:a_{1}x^{2}+a_{2}y^{2}+a_{3}z^{2}=r^{2},\text{ }a_{1},a_{2},a_{3}\in
\mathbb{R}^{+}\}.
\end{equation*}

The sectional curvature of the $S_{a_{1},a_{2},a_{3}}^{2}$ generated by the
non-degenerated plane $\{u,v\}$ is defined as%
\begin{equation}
K(x,y)=\frac{B(R(x,y)x,y)}{B(x,x)B(y,y)-B(x,y)^{2}},  \label{1}
\end{equation}%
where $R$ is the Riemannian curvature tensor of the $%
S_{a_{1},a_{2},a_{3}}^{2}$ given as%
\begin{equation}
R(X,Y)Z=\nabla _{X}\nabla _{Y}Z-\nabla _{Y}\nabla _{X}Z-\nabla _{\lbrack
X,Y]}Z.  \label{2}
\end{equation}%
The ellipsoid has the constant sectional curvature according to the
elliptical inner product. Therefore, the curvature tensor $R$ written as
follows%
\begin{equation}
R(X,Y)Z=B(Z,X)Y-B(Z,Y)X. \label{3}
\end{equation}%

Let $\beta $ be a unit speed curve on the $S_{a_{1},a_{2},a_{3}}^{2}$
defined by $\beta (s)=\varphi (\alpha (s))$. Then, the unit normal vector
field $Z$ along the surface $S_{a_{1},a_{2},a_{3}}^{2}$ defined by
\begin{equation*}
Z=\frac{\varphi _{u}\times _{E}\varphi _{v}}{\left\Vert \varphi _{u}\times
\varphi _{v}\right\Vert }.
\end{equation*}%
Since, $S_{a_{1},a_{2},a_{3}}^{2}$ is sphere according to the elliptical inner product, the
unit normal vector field $Z$ along the elliptical 2 sphere $S_{a_{1},a_{2},a_{3}}^{2}$ equal
to the position vector of the curve $\beta $. Then we found an orthonormal
frame $\{t=\beta ^{\prime },y=\beta \times _{E}\beta ^{\prime },\beta \}$
which is called the elliptical Darboux frame along the curve $\beta $.

\begin{proposition}
\label{Pro}Let $\beta $ be a unit speed curve on the elliptical 2 sphere $%
S_{a_{1},a_{2},a_{3}}^{2}$ with nowhere vanishing curvature. Then the corresponding Darboux
frame equations of $\beta $ is written as
\begin{eqnarray}
t^{\prime } &=&-\beta +k_{g_{E}}y,  \label{8} \\
\beta ^{\prime } &=&t,  \notag \\
y^{\prime } &=&-k_{g_{E}}t,  \notag
\end{eqnarray}%
where $k_{n_{E}}=-1$, $k_{g_{E}}=B(\beta ^{\prime \prime },y)$ and $\tau
_{r}=0$ are called the \emph{asymptotic curvature}, \emph{geodesic curvature}%
, and \emph{principal curvature }of $\beta $ on the $S_{a_{1},a_{2},a_{3}}^{2}$,
respectively.
\end{proposition}

\begin{proof}
Let $\beta $ be a curve on the elliptical 2 sphere $S_{a_{1},a_{2},a_{3}}^{2}$. First we
have the derivative of the position vector $\beta $ is equal to the tangent
vector. Secondly, we can write%
\begin{equation*}
t^{\prime }=\alpha t+\beta \beta +\delta y,
\end{equation*}%
where we found%
\begin{eqnarray*}
\alpha &=&B(t^{\prime },t)=0, \\
\beta &=&B(t^{\prime },\beta )=-1, \\
\delta &=&B(t^{\prime },y)=k_{g_{E}}.
\end{eqnarray*}%
Therefore we obtain
\begin{equation*}
t^{\prime }=-\beta +k_{g_{E}}y.
\end{equation*}%
Similar calculation above we calculate that%
\begin{equation*}
y^{\prime }=-k_{g_{E}}t.
\end{equation*}
\end{proof}

Since $\{t,\beta ,y\}$ is a right-handed orthonormal basis of $\mathbb{R}%
_{a_{1},a_{2},a_{3}}^{3}$, it is found as the following relation%
\begin{equation}
y\times _{E}t=\beta ,\text{ }\beta \times _{E}y=t,\text{ }\beta \times
_{E}t=-y.  \label{10}
\end{equation}

\begin{theorem}
\label{Teo}Suppose that $T$ is a skew symmetric matrix in the following form
\begin{equation*}
T=\Delta \left(
\begin{array}{ccc}
0 & -\frac{u_{3}}{a_{1}} & \frac{u_{2}}{a_{1}} \\
\frac{u_{3}}{a_{2}} & 0 & -\frac{u_{1}}{a_{2}} \\
-\frac{u_{2}}{a_{3}} & \frac{u_{1}}{a_{3}} & 0%
\end{array}%
\right),
\end{equation*}%
such that $u=(u_{1},u_{2},u_{3})\in $ $\mathbb{R}_{a_{1},a_{2},a_{3}}^{3}$
is a unit vector. Then, the matrix exponential%
\begin{equation*}
R_{\theta }^{B,u}=e^{T\theta }=I+\sin \theta T+(1-\cos \theta )T^{2},
\end{equation*}%
gives an elliptical rotation on the ellipsoid $%
a_{1}x^{2}+a_{2}y^{2}+a_{3}z^{2}=1,$ $a_{i}\in \mathbb{R}^{+}.$ Moreover,
the matrix $R_{\theta }^{B,u}$ can be expressed as:%
\begin{equation*}
\left(
\begin{array}{ccc}
a_{1}u_{1}^{2}+(1-a_{1}u_{1}^{2})\cos \theta & \frac{\Delta u_{3}\sin \theta
}{a_{1}}-a_{2}u_{1}u_{2}(\cos \theta -1) & \frac{\Delta u_{2}\sin \theta }{%
a_{1}}-a_{3}u_{1}u_{3}(\cos \theta -1) \\
\frac{\Delta u_{3}\sin \theta }{a_{2}}-a_{1}u_{1}u_{2}(\cos \theta -1) &
a_{2}u_{2}^{2}+(1-a_{2}u_{2}^{2})\cos \theta & -\frac{\Delta u_{1}\sin
\theta }{a_{2}}-a_{3}u_{2}u_{3}(\cos \theta -1) \\
-\frac{\Delta u_{2}\sin \theta }{a_{3}}-a_{1}u_{1}u_{3}(\cos \theta -1) &
\frac{\Delta u_{1}\sin \theta }{a_{3}}-a_{2}u_{2}u_{3}(\cos \theta -1) &
a_{3}u_{3}^{2}+(1-a_{3}u_{3}^{2})\cos \theta%
\end{array}%
\right),
\end{equation*}%
where $\Delta =\sqrt{a_{1}a_{2}a_{3}}$ and $\theta $ is the elliptical rotation angle and $u=(u_{1},u_{2},u_{3})$ is the rotation axis \cite{Ozd1}.
\end{theorem}

\section{Variations of Darboux frame for elliptical 2-sphere}

In this section, we examine the variational formulations of the curves on
the elliptical 2 sphere $S_{a_{1},a_{2},a_{3}}^{2}$ to describe the motion of a point on the
$S_{a_{1},a_{2},a_{3}}^{2}$. First, we study the features of variational vector field along
a curve and calculate the variational formulas for its Darboux curvatures.
After that, we give the connection between the Killing equations along a curve and the geometric variational formulas for curvatures of the curve on the
ellipsoid $S_{a_{1},a_{2},a_{3}}^{2}.$

\begin{lemma}
\label{lemma}Let $\varphi :U\subset \mathbb{E}^{2}\rightarrow \mathbb{E}%
^{3}, $ $\varphi (U)=S_{a_{1},a_{2},a_{3}}^{2}$ be an elliptical 2-sphere, $\beta
:I\subset \mathbb{R}\rightarrow U$ be a regular curve on $S_{a_{1},a_{2},a_{3}}^{2}$ and $K$ be a vector field along the curve $\beta $. Then, the variations of the
geodesic curvature function $k_{g_{E}}(s,w)$ and the speed function $v(s,w)$
at $w=0$ are calculated as follows:\newline
\begin{equation}
\begin{array}{l}
K(v)=\left. (\frac{\partial v}{\partial w}(s,w))\right\vert _{w=0}=-v\wp, \\
\\
K(k_{g_{E}})=\left. (\frac{\partial k_{g_{E}}}{\partial w}(s,w))\right\vert
_{w=0}=B(R(K,t)t+\nabla _{t}^{2}K,y)-\frac{1}{k_{g_{E}}}B(\beta ,R(K,t)t+\nabla
_{t}^{2}K)-2\wp(\frac{1}{k_{g_{E}}}+1).%
\end{array}
\label{16}
\end{equation}%
where $\wp =-B(\nabla _{t}K,t)$ and $R$ stands for the curvature tensor of $%
S_{a_{1},a_{2},a_{3}}^{2}.$
\end{lemma}

\begin{proof}
Let $\varphi :U\subset \mathbb{R}^{2}\rightarrow \mathbb{R}^{3},$ $\varphi
(U)=S_{a_{1},a_{2},a_{3}}^{2}$ be an elliptical 2-sphere and $\beta :I\subset \mathbb{R}%
\rightarrow U$ be a regular curve on the $S_{a_{1},a_{2},a_{3}}^{2}$. Provided that $K$ be a
vector field along the curve $\beta $ and the speed function satisfy $%
W=vt(s,w).$ If we use the covariant derivative, $\nabla _{K}$, of the speed
function $v(s,w)$ and if we use the elliptical inner product we obtain%
\begin{eqnarray}
K(v) &=&\left. (\frac{\partial v}{\partial w}(s,w))\right\vert _{w=0}
\label{17} \\
&=&K(B(W,W)^{1/2}) \\
&=&\frac{1}{v}B(\nabla _{K}W,W) \\
&=&vB(\nabla tK,t)  \notag \\
&=&v\wp ,  \notag
\end{eqnarray}%
where $\wp =-B(\nabla tK,t).$

On the other hand, from the elliptical Darboux frame equations we have
\begin{equation}
k_{g_{E}}=B(\nabla _{t}t,y).  \label{18}
\end{equation}%
If we use the covariant derivative $\nabla _{K}$ of the geodesic curvature $%
k_{g}(s,t),$ we calculate%
\begin{equation}
K(k_{g_{E}})=\left. (\frac{\partial k_{g_{E}}}{\partial w}(s,w))\right\vert
_{w=0}=B(\nabla _{K}\nabla _{t}t,y)+B(\nabla _{t}t,\nabla _{K}y).  \label{19}
\end{equation}%
Also, we have the following two equations%
\begin{equation}
\nabla _{K}\nabla _{t}t=R(K,t)t+\nabla _{t}\nabla _{K}t+\nabla _{\lbrack
K,t]}t  \label{20}
\end{equation}%
and%
\begin{equation}
\lbrack K,t]=-\frac{1}{v}K(v)t;\text{ }-\frac{1}{v}K(v)=\wp .  \label{23}
\end{equation}%
The eq.(\ref{20}), eq.(\ref{23}), and the equation $\nabla
_{K}t=[K,t]+\nabla _{t}K$ give the following equation%
\begin{equation}
\nabla _{t}\nabla _{K}t=\wp ^{\prime }t+\wp \nabla _{t}t+\nabla _{t}^{2}K.
\label{24}
\end{equation}%
Moreover, from the elliptical Darboux frame equations, we have
\begin{equation}
y=\frac{1}{k_{g_{E}}}(\nabla _{t}t+\beta ).  \label{a}
\end{equation}%
The covariant derivative $\nabla _{K}$ of eq.(\ref{a}) implies%
\begin{equation*}
\nabla _{K}y=\frac{-k_{g_{E}}^{\prime }}{k_{g_{E}}^{2}}\frac{\partial
k_{g_{E}}}{\partial w}(\nabla _{t}t+\beta )+\frac{1}{k_{g_{E}}}(\nabla
_{K}\nabla _{t}t+\nabla _{K}\beta );
\end{equation*}%
Since, $\frac{\partial k_{g_{E}}}{\partial w}=0$ the last equation read%
\begin{equation}
\nabla _{K}y=\frac{1}{k_{g_{E}}}(\nabla _{K}\nabla _{t}t+\nabla _{K}\beta ).
\label{b}
\end{equation}

Using the eqs. (\ref{19})-(\ref{b}) the expression for $K(k_{g})$ becomes%
\begin{equation*}
\begin{array}{l}
K(k_{g_{E}})=\left. (\frac{\partial k_{g_{E}}}{\partial t}(s,t))\right\vert
_{t=0}=B(R(K,t)t+\nabla _{t}^{2}K,y)-\frac{1}{k_{g_{E}}}B(\beta ,R(K,t)t+\nabla
_{t}^{2}K)-2\wp(\frac{1}{k_{g_{E}}}+1).%
\end{array}%
\end{equation*}
\end{proof}

\begin{proposition}
\label{Pro1} Suppose that $K(s)$ is the restriction to the curve $\beta (s)$ of a Killing vector field $K$ of $S_{a_{1},a_{2},a_{3}}^{2}$. Then the variations of the speed function and the elliptical geodesic curvature of $\beta $ satisfy:%
\begin{equation}
K(v)=K(k_{g_{E}})=0.  \label{30}
\end{equation}
\end{proposition}

\begin{proof}
It is clearly stated that the variations $K(v)$ and $K(k_{g_{E}})$ do not depend on only $K(s)$. Then the variation of $\beta (s)$ in the direction of $K(s)$
is given as follows:
\begin{equation}
\beta _{t}(s)=\beta (s,t)=:F_{L_{t}} (\beta (s)).  \label{31}
\end{equation}%
Then, the isometric nature of $F_{L_{t}} $ implies that $K(v)=K(k_{g_{E}})=0.$ Since the ellipsoid has a constant sectional curvature, the Proposition \ref{Pro1}  has a converse.
\end{proof}

\section{Magnetic trajectories of the point on the elliptical 2-sphere}

When a charged particle enters in a magnetic field, under the action of the magnetic field it traces a trajectory called as magnetic curve. A magnetic field $F$ on the $S_{a_{1},a_{2},a_{3}}^{2}$ is a closed 2-form. The
Lorentz force of the magnetic field $F$ related to the elliptical inner
product on $S_{a_{1},a_{2},a_{3}}^{2}$ is defined to be a skew-symmetric operator fulfilling
that

\begin{equation}
B(F_{L} (X),Y)=F(X,Y),
\end{equation}%
for all $X=(x_{1},x_{2},x_{3}),Y=(y_{1},y_{2},y_{3})\in \chi (S_{a_{1},a_{2},a_{3}}^{2}).$
Then, the mixed product of the vector fields $X,Y,Z\in \chi (S_{a_{1},a_{2},a_{3}}^{2})$ is
satisfy

\begin{equation}
B(X\times _{E}Y,Z)=dv_{g}(X,Y,Z),
\end{equation}%
where $dv_{g}$ denotes a volume on the $S_{a_{1},a_{2},a_{3}}^{2}$.

Since, divergence free vector fields and magnetic fields are in one-to-one
correspondence, the Lorentz force $F_{L} $ related to the magnetic field
$K$ can be given by the following formula%
\begin{equation}
F_{L} (X)=K\times _{E}X.  \label{33}
\end{equation}%
From these formulates the trajectories produced by the magnetic field $K$
which satisfy the following Lorentz force equation

\begin{equation}
F_{L} (t)=K\times _{E}t=\nabla _{t}t,  \label{34}
\end{equation}%
where $\nabla $ is the Levi-Civita connection.

The studies conducted in this field could be briefly summarized as follows:
First studies of the magnetic curves were carried out Riemannian surfaces
\cite{Bar2,Sun}, and then, \cite{Dru,Dru1} some authors have provided some
characterizations for these curves in three-dimensional spaces as $\mathbb{E}%
^{3}$ and $\mathbb{E}_{1}^{3}$. Moreover, variational approach on magnetic
flows of the Killing magnetic field in $3D$ has been developed. Furthermore,
Barros \textit{et al. }\cite{Bar1}, revealed that solution of the Lorentz
force equation were found to be Kirchhoff elastic rods. According to the
study of Cabrerizo \cite{Cab}, it has been presented that the magnetic flow
lines associated with Killing magnetic fields on the unit $2D$ and $3D$
sphere are found. Finally, Bozkurt(Ozdemir) \textit{et al}.\cite{Boz,Ozd}
described N-magnetic and B-magnetic curves as the trajectories of the
certain magnetic field and they have provided some specifications for
magnetic flows associated with Killing magnetic field in$3D$ Riemannian and
semi-Riemannian manifolds. In addition, variational vector field has been
used for the characterization of elastic curves. In a study \cite{Gur},
Gurbuz has stated Killing equations in order to characterize the elastic
curves.

In this context, we define the trajectories of the charged particles on the $%
S_{a_{1},a_{2},a_{3}}^{2}$. We derive the parametric representations of all magnetic curves
on the $S_{a_{1},a_{2},a_{3}}^{2}$ and present various examples in order to illustrate the
results of the work in terms of theory.

\begin{proposition}
\label{Pro2}Let $\beta :I\subset
\mathbb{R}\rightarrow U$ be a regular curve on $S_{a_{1},a_{2},a_{3}}^{2}$. Then, the
Lorentz force equations related to the elliptical Darboux frame $%
\{t,\beta ,y\}$ can be expressed as the following matrix form:%
\begin{equation}
\left[
\begin{array}{c}
F_{L} (t) \\
F_{L} (\beta ) \\
F_{L} (y)%
\end{array}%
\right] =\left[
\begin{array}{ccc}
0 & -1 & k_{g_{E}} \\
1 & 0 & \delta \\
-k_{g_{E}} & -\delta & 0%
\end{array}%
\right] \left[
\begin{array}{c}
t \\
\beta \\
y%
\end{array}%
\right],  \label{35}
\end{equation}%
where $k_{g_{E}}$ is geodesic curvature\emph{\ }of $\beta $.
\end{proposition}

\begin{proof}
Let $\beta :I\subset \mathbb{R}%
\rightarrow U$ be a regular curve on $S_{a_{1},a_{2},a_{3}}^{2}$. Then, from the eq.(\ref{34}), we
have%
\begin{equation}
F_{L} (t)=\nabla _{t}t=-\beta +k_{g_{E}}y.  \label{36}
\end{equation}%
Since $F_{L} (\beta )\in span\{t,\beta ,y\}$we can write%
\begin{equation}
F_{L} (\beta )=\lambda t+\mu \beta +\xi y . \label{37}
\end{equation}%
Then we calculate the functions $\lambda ,$ $\mu$, and $\xi $ satisfy%
\begin{eqnarray}
\lambda &=&B(F_{L} (\beta ),t)=0,  \label{38} \\
\mu &=&B(F_{L} (\beta ),\beta )=0,  \notag \\
\xi &=&B(F_{L} (\beta ),y)=\delta .  \notag
\end{eqnarray}%
As a consequence $F_{L} (\beta )$ found as%
\begin{equation}
F_{L} (\beta )=t+\delta y.  \label{39}
\end{equation}%
Finally, similar computations lead to%
\begin{equation}
F_{L} (y)=-k_{g_{E}}t-\delta \beta .  \label{40}
\end{equation}
\end{proof}

\begin{proposition}
\label{Pro3}Let $\beta :I\subset \mathbb{%
R}\rightarrow U$ be a regular curve on the elliptical 2-sphere $S_{a_{1},a_{2},a_{3}}^{2}$ and $K$ be a Killing
vector field along the curve $\beta$. Then $\beta $ is a magnetic
path related to the magnetic field $K$ on the condition that $K$ can be
expressed along $\beta $ as:%
\begin{equation}
K(s)=\delta t-k_{g_{E}}\beta -y,  \label{41}
\end{equation}%
where $k_{g_{E}}$ is geodesic curvature and the function $\delta (s)$
associated with each magnetic curve will be expressed as \emph{quasislope},
calculated obviously, with respect to the magnetic field $K$.
\end{proposition}

\begin{proof}
We may express the Killing magnetic vector field $K(s)$ in the following form
\begin{equation}
K(s)=at+b\beta +cy,  \label{42}
\end{equation}%
where $a$, $b$ and $c$ are certain functions along the magnetic curve and we
assume that $K$ does not vanish on $\beta $. If the eqs. (\ref{33})-(\ref%
{35}) are used, we calculate $a=\delta ,b=-k_{g_{E}}$ and $c=-1$.
Conversely, if we assume that eq.(\ref{42}) holds then we will obtain $%
K\times t=F_{L} (t)$. Thus, the curve $\beta $ is a magnetic trajectory related to the
magnetic vector field $K$.
\end{proof}

\begin{theorem}
\label{Teo1}Let $\varphi :U\subset \mathbb{E}^{2}\rightarrow \mathbb{E}^{3},$
$\varphi (U)=S_{a_{1},a_{2},a_{3}}^{2}$ be an elliptical 2-sphere, $\beta :I\subset \mathbb{%
R}\rightarrow U$ be a regular curve on $S_{a_{1},a_{2},a_{3}}^{2}$ and $K$ be a Killing
vector field along the curve $\beta $. Then $\beta $ is a magnetic
trajectory of the magnetic field $K$ on the condition that the elliptical
Darboux frame curvatures of the curve $\beta $ satisfy the following
differential equation:%
\begin{equation*}
k_{g_{E}}^{\prime \prime }+\delta k_{g_{E}}k_{g_{E}}^{\prime }=0,
\end{equation*}%
where $\delta $ is a constant along the magnetic trajectory $\beta$.
\end{theorem}

\begin{proof}
If $\beta $ is a magnetic trajectory of a magnetic field $K$ then $K$
satisfies eq.(\ref{41}). Taking the derivative of eq.(\ref{41}),
we found%
\begin{equation}
\nabla _{t}K=\delta ^{\prime }t+(-k_{g_{E}}^{\prime }-\delta )\beta
+(\delta k_{g_{E}})y.  \label{43}
\end{equation}%
Then, the equation, $K(v)=0$, in Lemma \ref{lemma} implies that $\delta $ is a
constant. If we take the derivation of eq.(\ref{43}), we obtain%
\begin{equation}
\nabla _{t}^{2}K=(-k_{g_{E}}^{\prime }-\delta -\delta k_{g_{E}}^{2})t+\delta
k_{g_{E}}^{\prime }y-k_{g}^{\prime \prime }\beta .  \label{45}
\end{equation}%
Moreover, we have found
\begin{equation}
R(K,t)t=(B(t,K)t-B(t,t)K).  \label{46}
\end{equation}%
Using the Darboux frame equations and eq.(\ref{41}) we deduce%
\begin{equation}
R(K,t)t=(k_{g_{E}}\beta +y).  \label{47}
\end{equation}%
Considering the eq.(\ref{45}) and eq.(\ref{47}) with the second equation in Lemma %
\ref{lemma} and the Proposition \ref{Pro1}, we have reached the following
equation%
\begin{equation}
k_{g_{E}}^{\prime \prime }+\delta k_{g_{E}}k_{g_{E}}^{\prime }=0.  \label{49}
\end{equation}%
The solution of the differential equation obtained as follows%
\begin{equation*}
k_{g_{E}}=c\text{ or }k_{g_{E}}=\sqrt{\frac{2c_{1}}{\delta }}\tanh (\sqrt{%
\frac{\delta c_{1}(s+c_{2})}{2}})\text{ for }\delta\neq 0, \text{ or }k_{g_{E}}=c_{3}s+c_{4}, \text{ for }\delta= 0.
\end{equation*}
\end{proof}

\begin{theorem}
\label{Teo2}Let $\beta $ be a magnetic curve on the elliptical 2-sphere
then the trajectory $\beta $ has one of the following representations:\newline
i. If $\beta $ has the constant curvature $k_{g}=c,$%
\begin{equation}
\beta (s)=\eta _{1}+\eta _{2}\sin (c^{2}+1)s+\eta _{3}\cos (c^{2}+1)s\text{%
, }k_{g}=c,  \label{56}
\end{equation}%
where $\eta _{1},\eta _{2},\eta _{3}\in \mathbb{R}_{a_{1},a_{2},a_{3}}^{3}$
and $c\in \mathbb{R}.$\newline
ii. If the geodesic curvature of the curve $\beta $ has the form $k_{g_{E}}=%
\sqrt{\frac{2c_{1}}{\delta }}\tanh (\sqrt{\frac{\delta c_{1}(s+c_{2})}{2}})$,%
\begin{equation}
\beta (s)=\left(
\begin{array}{l}
\mu _{1}hypergeom([\frac{1}{2},\frac{\frac{1}{2}I(\sqrt{2c_{1}}\mp 2\delta I)%
}{\delta },\frac{-\frac{1}{2}I(\sqrt{2c_{1}}\pm 2\delta I)}{\delta }],[\frac{%
1}{2}\frac{\pm 2\delta +\sqrt{-\delta -2c_{1}}}{\delta }, \\
-\frac{1}{2}\frac{\mp 2\delta +\sqrt{-\delta -2c_{1}}}{\delta }],\frac{1}{%
\cosh (\sqrt{\delta }(s+c_{2}))^{2}})\tanh (\sqrt{\delta }(s+c_{2})), \\
+\frac{\mu _{2}}{\sqrt{a_{2}}}\cosh (\sqrt{\delta }(s+c_{2}))^{\frac{\sqrt{%
-\delta -2c_{1}}-\delta )}{\delta }}\sinh (\sqrt{\delta }(s+c_{2})) \\
hypergeom([-\frac{1}{2}\frac{\mp \delta +\sqrt{-\delta -2c_{1}}}{\delta },-%
\frac{1}{2}\frac{-\sqrt{2\delta }I+\sqrt{-\delta -2c_{1}}\mp 2\delta }{%
\delta },-\frac{1}{2}\frac{\sqrt{2\delta }I+\sqrt{-\delta -2c_{1}}\mp
2\delta }{\delta }], \\
\lbrack -\frac{\mp \delta +\sqrt{-\delta -2c_{1}}}{\delta },-\frac{1}{2}%
\frac{\mp 2\delta +\sqrt{-\delta -2c_{1}}}{\delta }],\frac{1}{\cosh (\sqrt{%
\delta }(s+c_{2}))^{2}}) \\
+\frac{\mu _{3}}{\sqrt{a_{3}}}\cosh (\sqrt{\delta }(s+c_{2}))^{-\frac{\sqrt{%
-\delta -2c_{1}}+\delta )}{\delta }}\sinh (\sqrt{\delta }(s+c_{2})) \\
hypergeom([\frac{1}{2}\frac{\pm \delta +\sqrt{-\delta -2c_{1}}}{\delta },%
\frac{1}{2}\frac{\sqrt{2\delta }I+\sqrt{-\delta -2c_{1}}\pm 2\delta }{\delta
},\frac{1}{2}\frac{-\sqrt{2\delta }I+\sqrt{-\delta -2c_{1}}\pm 2\delta }{%
\delta }], \\
\lbrack \frac{\pm \delta +\sqrt{-\delta -2c_{1}}}{\delta },\frac{1}{2}\frac{%
\pm 2\delta +\sqrt{-\delta -2c_{1}}}{\delta }],\frac{1}{\cosh (\sqrt{\delta }%
(s+c_{2}))^{2}})%
\end{array}%
\right),  \label{61}
\end{equation}%
where $\mu _{1},\mu _{2},\mu _{3}\in \mathbb{R}_{a_{1},a_{2},a_{3}}^{3}$ and
$c_{1},c_{2}\in \mathbb{R}.$
\end{theorem}

\begin{proof}
Since $\beta $ is a curve on the elliptical 2-sphere $S_{a_{1},a_{2},a_{3}}^{2}$ the curve
has the following curvatures%
\begin{equation*}
k_{g_{E}},\text{ }k_{n}=1,\text{ }\tau _{g}=0.
\end{equation*}%
Then from the elliptical Darboux frame equation we have the following third
order differential equation%
\begin{equation*}
k_{g_{E}}\beta ^{\prime \prime \prime }-k_{g_{E}}^{\prime }\beta ^{\prime
\prime }+(k_{g_{E}}^{3}+k_{g_{E}})\beta ^{\prime }-k_{g_{E}}^{\prime
}\beta =0.
\end{equation*}%
The solution of the differential equation is presented in the equations (\ref%
{61}) and (\ref{61}).
\end{proof}

\begin{example}
Theorem 4.4 shows that we can find the following magnetic curve on the $%
S_{a_{1},a_{2},a_{3}}^{2}$ parameterized by%
\begin{equation*}
\beta (s)=(\frac{\cos 3s}{\sqrt{2a_{1}}},\frac{\sin 3s}{\sqrt{2a_{2}}},%
\frac{1}{\sqrt{2a_{3}}}),
\end{equation*}%
with choosing the $\eta _{1}=(0,0,\frac{1}{\sqrt{2a_{3}}}),$ $\eta _{2}=(0,%
\frac{1}{\sqrt{2a_{2}}},0)$, $\eta _{3}=(\frac{1}{\sqrt{2a_{1}}},0,0).$ Then
the elliptical curvature computed as%
\begin{equation*}
k_{g_{E}}=\sqrt{2}.
\end{equation*}%
Different combination of the $\eta _{i}$, $i=1,2,3$, gives different positions of the circular magnetic curves. The images of the magnetic curves on the ellipsoid $S_{2,8,4}^{2}$ are shown in Figure \ref{fig:1}.%%
%%%%%%%%%%%%%%%%%%%%%%%%%%%%%%%%%%%%%%%%%%%%%%%%%%%%%%%%%
 \begin{figure}[H]
 \centering 
% \subfigure[]{\includegraphics[width=0.30\textwidth]{Figure7.jpg}}
\includegraphics[width=0.29\textwidth]{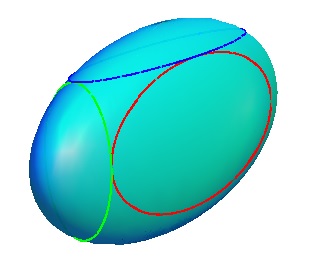}
 \caption{Charged particle motions along the curve $\beta$ \text{
on the }$S_{2,8,4}^{2}$\text{.} }\label{fig:1}
 \end{figure}
%%%%%%%%%%%%%%%%%%%%%%%%%%%%%%%%%%%%%%%%%%%%%%%%%%%%%%%%%%%%%
\end{example}

\begin{example}
Let $\beta $ be a curve on the elliptical 2-sphere and has the following
parametric representation
\begin{equation*}
\beta (s)=\left(
\begin{array}{c}
\frac{1}{\sqrt{a_{1}}}\cos (\int\limits_{0}^{t}\cos s^{2} ds)\sin (\int\limits_{0}^{t}\sin
s^{2} ds) \\
\frac{1}{\sqrt{a_{2}}}\sin (\int\limits_{0}^{t}\cos s^{2} ds)\sin
(\int\limits_{0}^{t}(\sin s^{2} ds) \\
\frac{1}{\sqrt{a_{3}}}\cos (\int\limits_{0}^{t}\sin s^{2} ds)%
\end{array}%
\right) .
\end{equation*}%
with the elliptical curvature%
\begin{equation*}
k_{g_{E}}=s.
\end{equation*}%
Then we can easily see that $\beta $ is a magnetic curve on the elliptical
2 sphere with $\delta =0$. In Figure \ref{fig:3}, we are visualized the image of the magnetic trajectory. %
%%%%%%%%%%%%%%%%%%%%%%%%%%%%%%%%%%%%%%%%%%%%%%%%%%%%%%%%%
 \begin{figure}[H]
 \centering 
% \subfigure[]{\includegraphics[width=0.30\textwidth]{Figure7.jpg}}
 \subfigure[]{\includegraphics[width=0.3\textwidth]{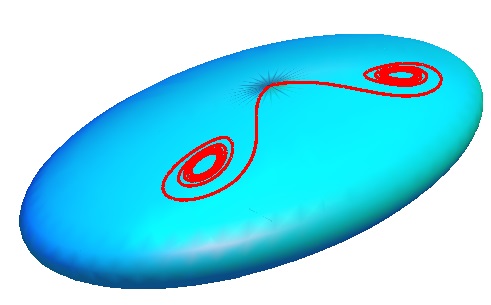}}
 \caption{\text{Charged particle motions along the curve }$\beta$ \text{ on the }$S_{1,2,4}^{2}$\text{,}}\label{fig:3}
 \end{figure}
%%%%%%%%%%%%%%%%%%%%%%%%%%%%%%%%%%%%%%%%%%%%%%%%%%%%%%%%%%%%%
\end{example}

\section{Helical trajectories of the point on the elliptical 2-sphere}
\subsection{Helical trajectories via variational methods}
In the following theorem, we introduce the helices on the elliptical
2-sphere defined as a curve whose tangent vector makes a constant angle with
a constant Killing vector field on the elliptical 2-sphere $S_{a_{1},a_{2},a_{3}}^{2}$.
\begin{theorem}
\label{Teo3}Let $\beta :I\subset \mathbb{R}\rightarrow U$ be a regular
curve on $S_{a_{1},a_{2},a_{3}}^{2}$ and $J$ be a Killing vector field along the curve $%
\beta $. Then $\beta $ is a helical trajectory of the magnetic field $K$
on the condition that the elliptical Darboux curvature of the curve $\beta $
has the form%
\begin{equation*}
k_{g_{E}}=\cot \theta,
\end{equation*}%
where $\theta $ satisfy%
\begin{equation}
\theta ^{\prime \prime }\sin ^{2}\theta -\omega \theta ^{\prime }\cos \theta
=0.
\end{equation}
\end{theorem}

\begin{proof}
If $\beta $ is a helical trajectory of the magnetic field $J$ then $J$
satisfies eq.(\ref{41}).%
\begin{equation}
J=\omega t+\cos \theta \beta +\sin \theta y.\label{x1}
\end{equation}%
Taking the derivative of eq.(\ref{41}), we found 
equations;%
\begin{eqnarray}
\nabla _{T}J &=&(\cos \theta -k_{g_{E}}\sin \theta )t+(-\theta ^{\prime
}\sin \theta -\omega )\beta +(\omega k_{g_{E}}+\theta ^{\prime }\cos \theta )y.
\end{eqnarray}%
Using the equation $J(v)=0$ in Lemma \ref{lemma} we present $\delta $ is a
constant.%
\begin{equation*}
k_{g_{E}}=\cot \theta.
\end{equation*}%
The differentiation of eq.(\ref{41}) is as follows;%
\begin{equation}
\begin{array}{l}
\nabla _{T}^{2}J=-\theta ^{\prime }\sin \theta -\omega -\omega
k_{g_{E}}^{2}-k_{g_{E}}\theta ^{\prime }\cos \theta )t \\
+(-\theta ^{\prime \prime }\sin \theta -\theta ^{\prime 2}\cos \theta
)\beta +(k_{g_{E}}^{\prime }\omega +\theta ^{\prime \prime }\cos \theta
-\theta ^{\prime 2}\sin \theta )y.%
\end{array}%
\end{equation}%
Moreover, we have found the following equation%
\begin{equation}
R(J,t)t=(B(t,J)t-B(t,t)J).
\end{equation}%
Using the Darboux frame equations and eq.(\ref{41}) we deduce%
\begin{equation}
R(J,t)t=(k_{g_{E}}\beta +y).
\end{equation}%
Considering the eq.(\ref{45}) and eq.(\ref{47}) with the second equation in Lemma \ref{lemma} and the Proposition \ref{Pro1}, we have reached the following
equation%
\begin{equation}
\theta ^{\prime \prime }\sin ^{2}\theta -\omega \theta ^{\prime }\cos \theta
=0.
\end{equation}
\end{proof}

\subsection{Helical trajectories via elliptical rotations}

In this subsection, we proved that helices on the $S_{a_{1},a_{2},a_{3}}^{2}$ can be described
by the motion of a point on a great ellipse of the $S_{a_{1},a_{2},a_{3}}^{2}$ rolling
without slipping on a fixed ellipse of the $S_{a_{1},a_{2},a_{3}}^{2}$.
\begin{theorem}
Let $\beta$ be generalized magnetic helix in $S_{a_{1},a_{2},a_{3}}^{2}$. Then $\beta$ can be expressed in terms of the elliptical motion formula as follows:
\begin{equation*}
   \beta(t)=R_{t }^{B,u}R_{\alpha }^{B,u}P,
\end{equation*}
where $R_{t }^{B,u}$, $R_{\alpha }^{B,u}$, $\alpha =\arccos k$, are  elliptical rotation matrices and $P$ is a point on the rolling great ellipse of the $S_{a_{1},a_{2},a_{3}}^{2}$.
\end{theorem}
\begin{proof}
Let $P$ be a point on the rolling ellipse centered at the origin in the $xoy-$%
plane given by

\begin{equation*}
P=(\frac{\cos kt}{\sqrt{a_{1}}},\frac{\sin kt}{\sqrt{a_{2}}},0).
\end{equation*}%
From the Theorem \ref{Teo} may take the axis as $u=(0,\frac{1}{\sqrt{a_{2}}}%
,0).$ Then the rotation matrix around the $y-axis$ with the angle $\alpha $
can be expressed as:%
\begin{equation*}
R_{\alpha }^{B,u}=\left(
\begin{array}{ccc}
\cos \alpha  & 0 & \frac{\sqrt{a_{3}}}{\sqrt{a_{1}}}\sin \alpha  \\
0 & 1 & 0 \\
-\frac{\sqrt{a_{1}}}{\sqrt{a_{3}}}\sin \alpha  & 0 & \cos \alpha
\end{array}%
\right).
\end{equation*}%
If we take the axis as $u=(0,0,\frac{1}{\sqrt{a_{3}}})$ in Theorem \ref{Teo}
then we express the rotation matrix around the $z-axis$ with the angle $t$
as:
\begin{equation*}
R_{t}^{B,u}=\left(
\begin{array}{ccc}
\cos t & \frac{\sqrt{a_{2}}}{\sqrt{a_{1}}}\sin t & 0 \\
-\frac{\sqrt{a_{1}}}{\sqrt{a_{2}}}\sin t & \cos t & 0 \\
0 & 0 & 1%
\end{array}%
\right).
\end{equation*}%
First we rotate this ellipse around the $y-axis$ with the angle $\alpha
=\arccos k$ we obtain:%
\begin{eqnarray*}
&&\left(
\begin{array}{ccc}
\cos \alpha  & 0 & \frac{\sqrt{a_{3}}}{\sqrt{a_{1}}}\sin \alpha  \\
0 & 1 & 0 \\
-\frac{\sqrt{a_{1}}}{\sqrt{a_{3}}}\sin \alpha  & 0 & \cos \alpha
\end{array}%
\right) \left(
\begin{array}{c}
\frac{\cos kt}{\sqrt{a_{1}}} \\
\frac{\sin kt}{\sqrt{a_{2}}} \\
0%
\end{array}%
\right)  
=(\frac{k\cos kt}{\sqrt{a_{1}}},\frac{\sin kt}{\sqrt{a_{2}}},\frac{\sqrt{%
1-k^{2}}\cos kt}{\sqrt{a_{3}}}).
\end{eqnarray*}%
Then, if we rotate this ellipse by an angle $t$ around the $z-axis$ we have
\begin{equation*}
\beta(t) =\left(
\begin{array}{ccc}
\cos t & \frac{\sqrt{a_{2}}}{\sqrt{a_{1}}}\sin t & 0 \\
-\frac{\sqrt{a_{1}}}{\sqrt{a_{2}}}\sin t & \cos t & 0 \\
0 & 0 & 1%
\end{array}%
\right) \left(
\begin{array}{c}
\frac{k\cos kt}{\sqrt{a_{1}}} \\
\frac{\sin kt}{\sqrt{a_{2}}} \\
\frac{\sqrt{1-k^{2}}\cos kt}{\sqrt{a_{3}}}%
\end{array}%
\right) .
\end{equation*}
\end{proof}
Consequently, the parametric equations of the helices on the elliptical
2-sphere has the following parameterizations%
\begin{equation*}
\beta(t)=R\left(
\begin{array}{c}
\frac{k\cos t\cos kt+\sin t\sin kt}{\sqrt{a_{1}}} \\
\frac{k\sin t\cos kt-\cos t\sin kt}{\sqrt{a_{2}}} \\
\frac{\sqrt{1-k^{2}}\cos kt}{\sqrt{a_{3}}}%
\end{array}%
\right)
\end{equation*}%
or the other parameterizations
\begin{equation}
\beta(t)=R\left(
\begin{array}{c}
\frac{1}{2}\frac{(1+k)\cos (1-k)t-(1-k)\cos (1+k)t}{\sqrt{a_{1}}} \\
\frac{1}{2}\frac{(1+k)\sin (1-k)t-(1-k)\sin (1+k)t}{\sqrt{a_{2}}} \\
\frac{\sqrt{1-k^{2}}\cos kt}{\sqrt{a_{3}}}%
\end{array}%
\right).  \label{y}
\end{equation}

The elliptical curvature of the helix $\beta$ calculated as
\begin{equation}
k_{g_{E}}(s)=\cot ks.  \label{x}
\end{equation}%
The curvature in eq.(\ref{x}) satisfy the eq.(\ref{49}) in the Theorem \ref%
{Teo3}. This proves that the eq.(\ref{y}) is a generalized magnetic helix in the elliptical 2 sphere $S_{a_{1},a_{2},a_{3}}^{2}$. 

\begin{corollary}
If the general helices $\beta $ has the curvature $k_{g_{E}}=\cot ks$ and
the quassislope $\delta =2k$ then the curve $\beta $ is a magnetic curve on
the elliptical 2-sphere.
\end{corollary}
\begin{proof}
Let $\beta $ be a curve with the curvature $k_{g_{E}}=\cot ks.$ on the
other hand, if $\beta $ is a magnetic curve on the elliptical 2-sphere $%
S_{a_{1},a_{2},a_{3}}^{2}$ then $\beta $ satisfy the eq.(\ref{49}). It is easily stated that $%
k_{g_{E}}$ satisfy the eq.(\ref{49}) for $\delta =2k$.
\end{proof}
\begin{theorem}
If the general helices $\beta $ has the curvature $k_{g_{E}}=\cot \theta$,  $\theta=as+b$, then the curve $\beta $ is a magnetic curve on
the elliptical 2-sphere with the quassislope $\delta =2a$ .
\end{theorem}
\begin{example}
The images of the
helices corresponding to the different value of the $k$ illustrated in Figure \ref{fig:4} and Figure \ref{fig:5}.
%%%%%%%%%%%%%%%%%%%%%%%%%%%%%%%%%%%%%%%%%%%%%%%%%%%%%%%%%
 \begin{figure}[H]
 \centering 
 \subfigure[]{\includegraphics[width=0.2\textwidth]{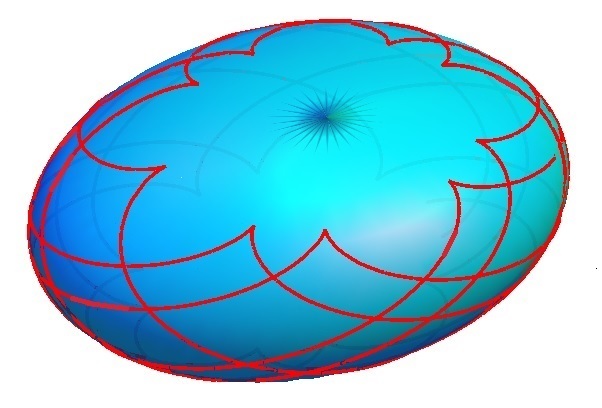}}
 \subfigure[]{\includegraphics[width=0.17\textwidth]{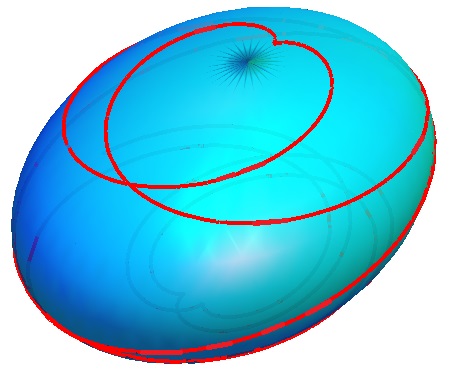}}
 \subfigure[]{\includegraphics[width=0.2\textwidth]{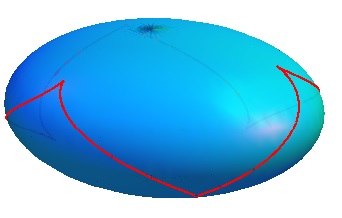}}
 \subfigure[]{\includegraphics[width=0.2\textwidth]{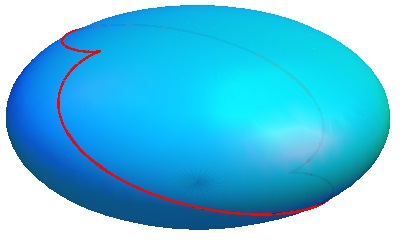}}
 \caption{Helical trajectories for, $k=0.56$,\text{ }$k=0.17$\text{, }$k=0.75$,\text{ }$k=0.5$, respectively.}\label{fig:4}
 \end{figure}
%%%%%%%%%%%%%%%%%%%%%%%%%%%%%%%%%%%%%%%%%%%%%%%%%%%%%%%%%%%%%
%%%%%%%%%%%%%%%%%%%%%%%%%%%%%%%%%%%%%%%%%%%%%%%%%%%%%%%%%
 \begin{figure}[H]
 \centering 
 \subfigure[]{\includegraphics[width=0.2\textwidth]{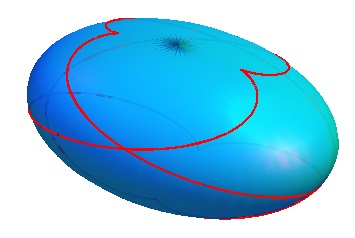}}
 \subfigure[]{\includegraphics[width=0.2\textwidth]{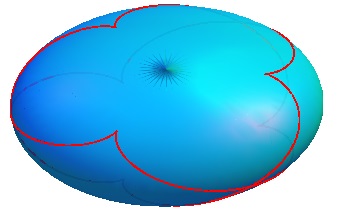}}
 \subfigure[]{\includegraphics[width=0.2\textwidth]{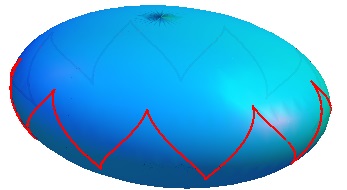}}
 \subfigure[]{\includegraphics[width=0.2\textwidth]{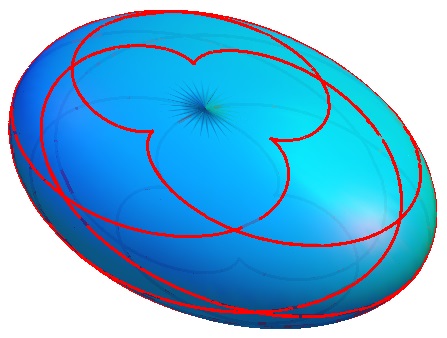}}
 \caption{Helical trajectories for, $k=0.4$,\text{ }$k=0.6$, $k=0.9$,\text{ }$k=0.3$, respectively.}\label{fig:5}
 \end{figure}
%%%%%%%%%%%%%%%%%%%%%%%%%%%%%%%%%%%%%%%%%%%%%%%%%%%%%%%%%%%%%

\end{example}
In the following subsection we define the satellite curves on the elliptical
2-sphere by using the elliptical rotation matrix in the Theorem \ref{Teo}.

\subsection{Elliptical satellite curves on the ellipsoid}

trajectories of a point $P$ on a given great ellipse of the elliptical
2-sphere rotating around the its axes, this axis being itself in uniform
rotation around an axis passing by the center of the ellipse.

The satellite curves include the helices on the elliptical 2-sphere when $%
k=-\cos \alpha $; it is the case where the curve has cuspidal points. Now we
obtain the parametric representation of the satellite curve according to the
elliptical motion. Let we take an arbitrary po\i int on the rolling ellipse
centered at $(0,0,0)$ in the $xoy-$plane given by%
\begin{equation*}
P=(\frac{\cos kt}{\sqrt{a_{1}}},\frac{\sin kt}{\sqrt{a_{2}}},0)
\end{equation*}%
First, if we rotate the point $P$ around the $y-axis$ we obtain:%
\begin{eqnarray*}
&&\left(
\begin{array}{ccc}
\cos \alpha & 0 & \frac{\sqrt{a_{3}}}{\sqrt{a_{1}}}\sin \alpha \\
0 & 1 & 0 \\
-\frac{\sqrt{a_{1}}}{\sqrt{a_{3}}}\sin \alpha & 0 & \cos \alpha
\end{array}
\right) \left(
\begin{array}{c}
\frac{\cos kt}{\sqrt{a_{1}}} \\
\frac{\sin kt}{\sqrt{a_{2}}} \\
0
\end{array}
\right) \\
&=&(\frac{\cos \alpha \cos kt}{\sqrt{a_{1}}},\frac{\sin kt}{\sqrt{a_{2}}},
\frac{\sin \alpha \cos kt}{\sqrt{a_{3}}}).
\end{eqnarray*}
Then, if we rotate the point $(\frac{\cos \alpha \cos kt}{\sqrt{a_{1}}},
\frac{\sin kt}{\sqrt{a_{2}}},\frac{\sin \alpha \cos kt}{\sqrt{a_{3}}})$ by
an angle $t$ around the $z-axis$, we have
\begin{equation*}
\beta(t) =\left(
\begin{array}{ccc}
\cos t & -\frac{\sqrt{a_{2}}}{\sqrt{a_{1}}}\sin t & 0 \\
\frac{\sqrt{a_{1}}}{\sqrt{a_{2}}}\sin t & \cos t & 0 \\
0 & 0 & 1
\end{array}
\right) \left(
\begin{array}{c}
\frac{\cos \alpha \cos kt}{\sqrt{a_{1}}} \\
\frac{\sin kt}{\sqrt{a_{2}}} \\
\frac{\sin \alpha \cos kt}{\sqrt{a_{3}}}
\end{array}
\right) .
\end{equation*}

Finally, the parametric equations of the satellite curve on the elliptical
2-sphere calculated as follows:
\begin{equation*}
\beta(t) =R\left(
\begin{array}{c}
\frac{\cos t\cos kt\cos \alpha -\sin t\sin kt}{\sqrt{a_{1}}} \\
\frac{\sin t\cos kt\cos \alpha +\cos t\sin kt}{\sqrt{a_{2}}} \\
\frac{\cos kt\sin \alpha }{\sqrt{a_{3}}}
\end{array}
\right) .
\end{equation*}
In the following figures we give various satellite curves images on the
elliptical 2-sphere.
%%%%%%%%%%%%%%%%%%%%%%%%%%%%%%%%%%%%%%%%%%%%%%%%%%%%%%%
 \begin{figure}[H]
 \centering 
 \subfigure[]{\includegraphics[width=0.2\textwidth]{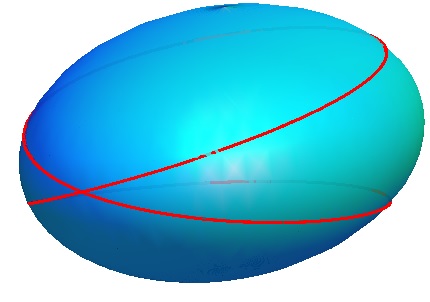}}
 \subfigure[]{\includegraphics[width=0.2\textwidth]{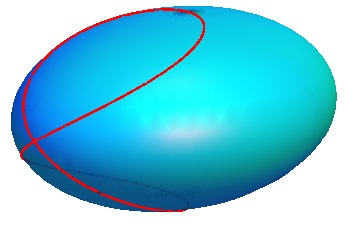}}
 \subfigure[]{\includegraphics[width=0.2\textwidth]{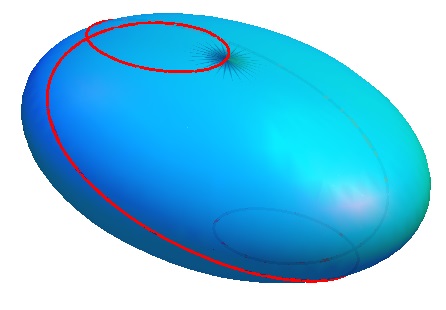}}
 \subfigure[]{\includegraphics[width=0.2\textwidth]{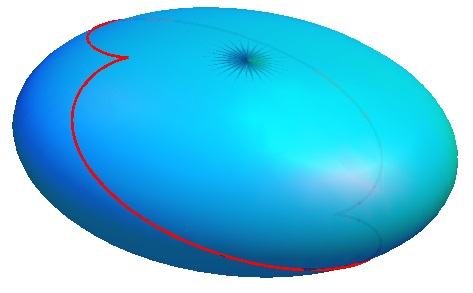}}
 \caption{$\alpha =\frac{1}{\sqrt{2}}$,\text{ }$k=1$;\text{ }$\alpha =\frac{\pi }{2}$,\text{
}$k=1$;\text{ }$\alpha =\frac{\pi }{2}$,\text{ }$k=\frac{1}{2}$;\text{ }$\alpha
=2.1$,\text{ } $k=\frac{1}{2}$, respectively. }\label{fig:6}
 \end{figure}
%%%%%%%%%%%%%%%%%%%%%%%%%%%%%%%%%%%%%%%%%%%%%%%%%%%%%%%%%%

%%%%%%%%%%%%%%%%%%%%%%%%%%%%%%%%%%%%%%%%%%%%%%%%%%%%%%%%%
 \begin{figure}[H]
 \centering 
 \subfigure[]{\includegraphics[width=0.2\textwidth]{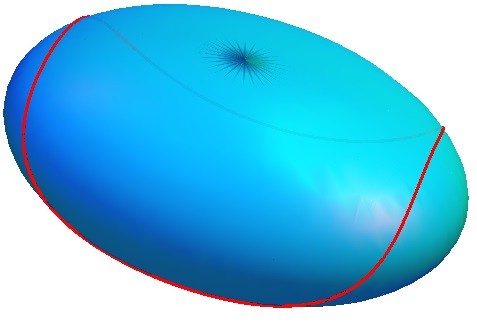}}
 \subfigure[]{\includegraphics[width=0.2\textwidth]{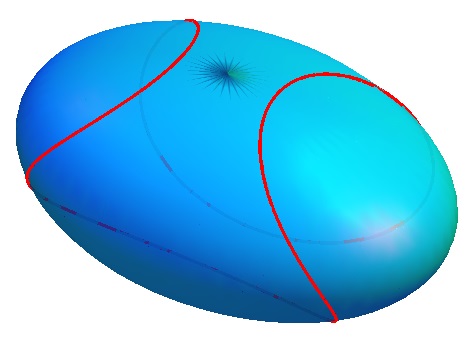}}
 \subfigure[]{\includegraphics[width=0.2\textwidth]{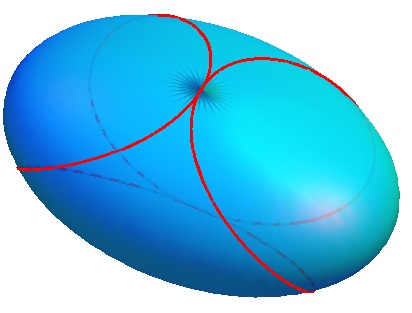}}
 \subfigure[]{\includegraphics[width=0.2\textwidth]{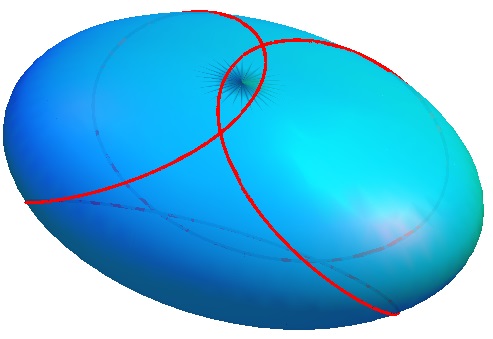}}
 \caption{$\alpha =2.5,\text{ }k=2;\text{ }\alpha =1.8,\text{ }k=2;\text{ }\alpha =1.57,
\text{ }k=2;\text{ }\alpha =1.5,\text{ }k=2$, respectively.}\label{fig:7}
 \end{figure}
%%%%%%%%%%%%%%%%%%%%%%%%%%%%%%%%%%%%%%%%%%%%%%%%%%%%%%%%%%%
In the following subsection we obtain the parametric representation of the
cycloid on the elliptical 2-sphere using the elliptical motion.

\subsection{Elliptical cycloid on the ellipsoid}

A cycloid on the elliptical 2 sphere is the trajectories of a point on an
ellipse rolling without slipping on a fixed ellipse on the elliptical 2
sphere such that the angle between the two ellipse remaining constant.
Suppose that the angle between the two ellipse equal to $\omega $ and the
fixed ellipse lies the plane.$xoy$. Since the rolling circle rolls without
sliding, its angular displacement $t$ around its center is $\frac{a}{b}t,$
here $a$ is the radius of the fixed ellipse, $b=\frac{a}{q}$ that of the
moving ellipse. A point on a copy of the rolling ellipse centered at $
(0,0,0) $ at in the $xoy-$plane given by
\begin{equation*}
P=(\frac{b}{\sqrt{a_{1}}}(1-\cos qt),\frac{b}{\sqrt{a_{2}}}(\sin qt),0)
\end{equation*}
First rotate this ellipse around the $y-axis$ we obtain:
\begin{eqnarray*}
&&\left(
\begin{array}{ccc}
\cos \omega & 0 & \frac{\sqrt{a_{3}}}{\sqrt{a_{1}}}\sin \omega \\
0 & 1 & 0 \\
-\frac{\sqrt{a_{1}}}{\sqrt{a_{3}}}\sin \omega & 0 & \cos \omega
\end{array}
\right) \left(
\begin{array}{c}
\frac{b}{\sqrt{a_{1}}}(1-\cos qt) \\
\frac{b}{\sqrt{a_{2}}}(\sin qt) \\
0
\end{array}
\right) \\
&=&(\frac{b}{\sqrt{a_{1}}}(1-\cos qt)\cos \omega ,\frac{b}{\sqrt{a_{2}}}
(\sin qt),-\frac{b}{\sqrt{a_{1}}}(1-\cos qt)\sin \omega ).
\end{eqnarray*}
Now translate the ellipse over a distance $a$ along the $x-axis$ to get:
\begin{equation*}
(\frac{a}{\sqrt{a_{1}}}+\frac{b}{\sqrt{a_{1}}}(1-\cos qt)\cos \omega ,\frac{b
}{\sqrt{a_{2}}}(\sin qt),-\frac{b}{\sqrt{a_{1}}}(1-\cos qt)\sin \omega ).
\end{equation*}
Finally, rotate this ellipse by an angle $t$ around the $z-axis$:
\begin{equation*}
\beta(t) =\left(
\begin{array}{ccc}
\cos t & \frac{\sqrt{a_{2}}}{\sqrt{a_{1}}}\sin t & 0 \\
-\frac{\sqrt{a_{1}}}{\sqrt{a_{2}}}\sin t & \cos t & 0 \\
0 & 0 & 1
\end{array}
\right) \left(
\begin{array}{c}
\frac{a}{\sqrt{a_{1}}}+\frac{b}{\sqrt{a_{1}}}(1-\cos qt)\cos \omega \\
\frac{b}{\sqrt{a_{2}}}(\sin qt) \\
-\frac{b}{\sqrt{a_{1}}}(1-\cos qt)\sin \omega
\end{array}
\right)
\end{equation*}

Then we obtain the parametric equations of the cycloid as
\begin{equation*}
\beta(t) =\left(
\begin{array}{c}
(\frac{a}{\sqrt{a_{1}}}+\frac{b}{\sqrt{a_{1}}}(1-\cos qt)\cos \omega )\cos t+
\frac{b}{\sqrt{a_{2}}}(\sin qt)\sin t \\
(\frac{a}{\sqrt{a_{1}}}+\frac{a}{\sqrt{a_{1}}}(1-\cos qt)\cos \omega )\sin t-
\frac{b}{\sqrt{a_{2}}}(\sin qt)\cos t \\
-\frac{b}{\sqrt{a_{1}}}(1-\cos qt)\sin \omega
\end{array}
\right) .
\end{equation*}
When $\omega =\pi $, we get the hypocycloid, and when $\omega =0$, the
epicycloid;
%%%%%%%%%%%%%%%%%%%%%%%%%%%%%%%%%%%%%%%%%%%%%%%%%%%%%
 \begin{figure}[H]
 \centering 
 \subfigure[]{\includegraphics[width=0.25\textwidth]{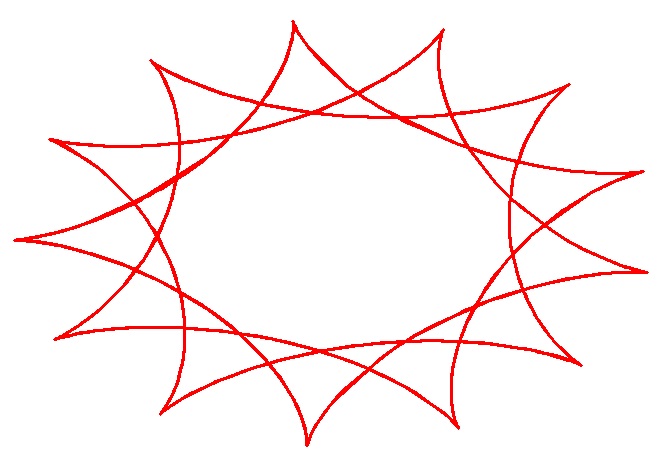}}
 \subfigure[]{\includegraphics[width=0.25\textwidth]{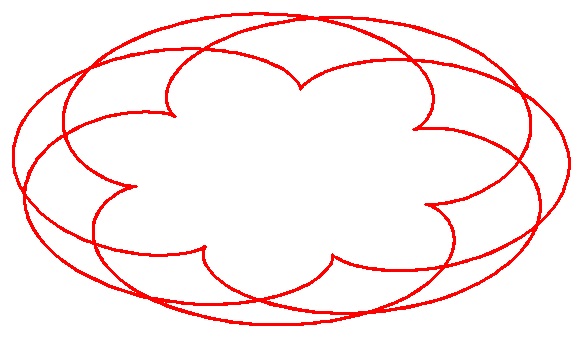}}
 \caption{hypocycloid and epicycloid, respectively. }\label{fig:8}
 \end{figure}
%%%%%%%%%%%%%%%%%%%%%%%%%%%%%%%%%%%%%%%%%%%%%%%%%%%%%%%%%
apart from these two cases, the cycloid is traced on the sphere
corresponding to both the base and the rolling ellipses, hence its name of
spherical cycloid. If $\omega =\arccos (-\frac{b}{a})$ the cycloid is
spherical cycloid with radius $a.$
%%%%%%%%%%%%%%%%%%%%%%%%%%%%%%%%%%%%%%%%%%%%%%%%%%%%%%%
 \begin{figure}[H]
 \centering 
 \subfigure[]{\includegraphics[width=0.2\textwidth]{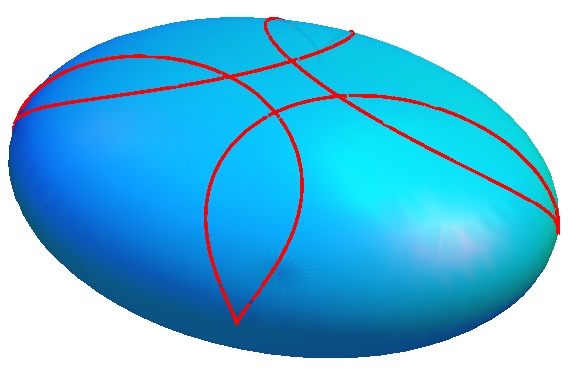}}
 \subfigure[]{\includegraphics[width=0.2\textwidth]{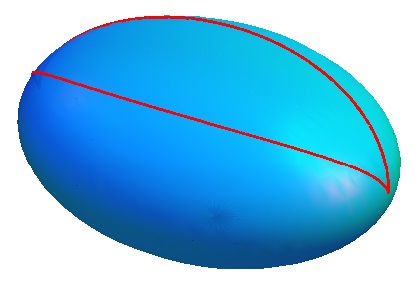}}
 \subfigure[]{\includegraphics[width=0.2\textwidth]{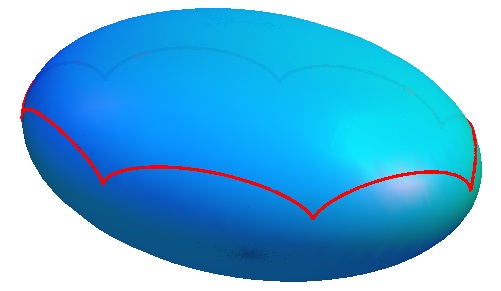}}
 \caption{$a=4,b=3;\text{ }a=4,b=2;\text{ }a=7,b=1$, respectively.}\label{fig:9}
 \end{figure}
%%%%%%%%%%%%%%%%%%%%%%%%%%%%%%%%%%%%%%%%%%%%%%%%%%%%%%%%%%
%%%%%%%%%%%%%%%%%%%%%%%%%%%%%%%%%%%%%%%%%%%%%%%%%%%%%%%%%
 \begin{figure}[H]
 \centering 
 \subfigure[]{\includegraphics[width=0.2\textwidth]{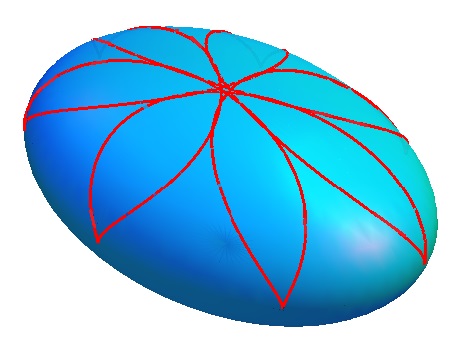}}
 \subfigure[]{\includegraphics[width=0.2\textwidth]{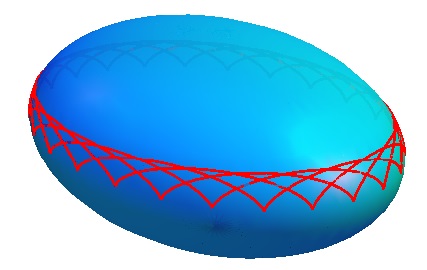}}
 \subfigure[]{\includegraphics[width=0.2\textwidth]{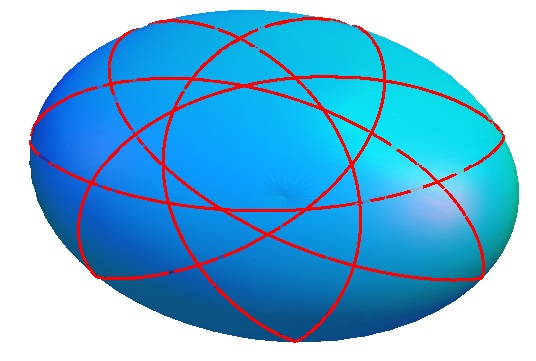}}
 \subfigure[]{\includegraphics[width=0.2\textwidth]{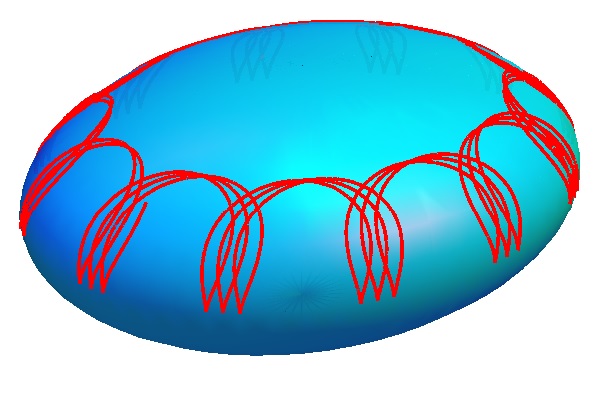}}
 \caption{$a=7,b=5;\text{ }a=25,b=4;\text{ }a=7,b=3, \text{ }a=1.09,b=1$, respectively.}\label{fig:10.}
 \end{figure}

\section{Conclusions}
In the present paper, the magnetic trajectories of a charged particle are investigated. To determine these trajectories the Killing equations are computed in terms of the Darboux frame along the ellipsoidal curve. Killing equations make it possible to interpret the motion of the charged particle on the elliptical 2-sphere $S_{a_{1}, a_{2}, a_{3}}^{2}$. By means of this approach, we determine the parametric representations of the magnetic trajectories. Besides, various examples of these trajectories are illustrated. Furthermore, if $ a_{1} = a_{2} = a_{3} = 1 $, the given characterizations are obtained on Euclidean 2-sphere $S^{2}$. This type of magnetic trajectories on the spheres was obtained for the first time with this study.

Next, we give some geometrical and physical interpretations of magnetic trajectories on the elliptical 2-sphere.
A charged particle which moves along a curve in a magnetic field $K$ experiences a force. This force is called as the Lorentz force, and satisfies that the equation
\begin{equation}
    F_{L}(t) = K\times t,
\end{equation}

where $K$ is the magnetic field and $T$ is velocity vector of the charged particle. 

In a static magnetic field, the trajectories are helical lines and the axis of the helix is parallel to the magnetic vector field. Since the Lorentz force is always perpendicular to the motion of the charged particle, if the charged particle moves parallel to magnetic field then the Lorentz force acts zero. If the charged particle moves perpendicular to magnetic field then the Lorentz force at the largest value and the particle traces a circular path. If the charged particle moves in a magnetic field with making a constant angle(not perpendicular) to magnetic field then experiences the Lorentz force and the particle traces a helical path in Euclidean 3-space (see, \cite{Bar,Bar1,Bar2}). 

In view of our approach we obtain the following characterizations:

i. If the velocity vector $t$ is parallel to the magnetic field $K$, then the particle moves parallel to the magnetic field since the Lorentz force will be zero. 

ii. If the velocity vector $t$ is perpendicular to the magnetic field $K$, the Lorentz force takes its maximum value and follows a trajectory that the elliptical geodesic curvature is a linear function. In this case, for example, the charged particle can be moves along the elliptical trajectory or spiral trajectory which related to its energy. Namely, first the particle enters to the magnetic field perpendicularly and follows an elliptical trajectory, and as the energy decreases, then this trajectory turns into a spiral trajectory or the particle first begins to follow a spiral trajectory, and as the energy decreases, then it can continue its motion by following an elliptical trajectory with constant energy. 
%In Figure \ref{fig:10}, we are visualised the charged particle motion for the constant angle $\delta=\pi/2$.
%%%%%%%%%%%%%%%%%%%%%%%%%%%%%%%%%%%%%%%%%%%%%%%%%%%%%%%%%
 %\begin{figure}[H]
 %\centering 
% \includegraphics[width=0.5\textwidth]{Figure31.jpg}
 %\caption{Charged particle motion along the Euler spiral on the %$S_{4,2,9}^{2}$. }\label{fig:10}
 %\end{figure}
%%%%%%%%%%%%%%%%%%%%%%%%%%%%%%%%%%%%%%%%%%%%%%%%%%%%%%%%%%%%%

iii. If the particle enters the magnetic field $K$ at a fixed angle(not perpendicular), the particle follows a trajectory that has the elliptical curvature $k_{g_{E}}=\sqrt{\frac{2c_{1}}{\delta }}\tanh (\sqrt{%
\frac{\delta c_{1}(s+c_{2})}{2}})$ or $k_{g_{E}}=cot\theta$, $\theta$ is a linear function, under the influence of the Lorentz force. For $k_{g_{E}}=cot\theta$, $\theta$ is linear, the charged particle traces a helical path on the ellipsoid. This helix makes a constant angle with the following Killing magnetic fields  
\begin{equation}
J=\omega t+\cos \theta \beta +\sin \theta y \textit{ and } V(s)=\delta t-k_{g_{E}}\beta -y,
\end{equation}%
Namely, these helices are curves that make a constant slope with respect to the plane $span\{V,J\}$, traced on a sphere. 
%In Figure \ref{fig:11}, we are visualised the charged particle motion for the constant angle $\delta=2k$, $k\neq0$.

%%%%%%%%%%%%%%%%%%%%%%%%%%%%%%%%%%%%%%%%%%%%%%%%%%%%%%%%%
% \begin{figure}[H]
 %\centering 
 %\includegraphics[width=0.42\textwidth]{Figure30.jpg}
 %\caption{Charged particle motion along the helical path on the %$S_{4,2,9}^{2}$. }\label{fig:11}
% \end{figure}
%%%%%%%%%%%%%%%%%%%%%%%%%%%%%%%%%%%%%%%%%%%%%%%%%%%%%%%%%%%%%

\textbf{Acknowledgement}: This research has been supported by Amasya University Scientific Research Projects Coordination Unit. Project Number: FMB-BAP 20-0452.

%\textbf{data availability statement:} Data sharing is not applicable to this article as no datasets were generated or analyzed during the current study.

%textbf{Authors’ Contributions:} Z.Ö. and F.A. equally contributed to the design and implementation of the research, to the analysis of the results, and to the writing of the manuscript.

%\textbf{Compliance with Ethical Standards:}
%Conflict of interest: The authors declare no conflicts of interest.

\end{document}